\documentclass{amsart}
\usepackage{amssymb,amsmath,amscd, amsfonts,latexsym}
\usepackage{graphicx}
\usepackage{amssymb}
\usepackage{epstopdf}

\usepackage{verbatim}  

\usepackage[all]{xy}
\usepackage{color}
\newtheorem{theorem}{Theorem}

\newtheorem{lemma}{Lemma}
\newtheorem{proposition}{Proposition}
\theoremstyle{remark}
\newtheorem{remark}{Remark}

\newcommand{\bbQ}{\mathbb Q}

\newcommand{\bbZ}{\mathbb Z}

\title[\title{On a variant of Pillai's problem}]{On a variant of Pillai's problem}

\author[K. C. Chim]{Kwok Chi Chim}

\author[I. Pink]{Istv\'an Pink}

\author[V. Ziegler]{Volker Ziegler}

\thanks{The first author was supported by the Austrian Science Fund (FWF) under the projects P24574, P26114 and W1230.
The second and the third authors were supported by the Austrian Science Fund (FWF) under the project P 24801-N26. }



\address{K. C. Chim \newline
         \indent Institute of Analysis and Number Theory, Graz University of Technology \newline
         \indent Steyrergasse 30/II \newline
         \indent A-8010 Graz, Austria}

\email{chim\char'100math.tugraz.at}

\address{I. Pink \newline
         \indent Institute of Mathematics, University of Debrecen \newline
         \indent H-4010 Debrecen, P.O. Box 12, Hungary \newline
         \indent and \newline
         \indent University of Salzburg \newline
         \indent Hellbrunnerstrasse 34/I \newline
         \indent A-5020 Salzburg, Austria}

\email{pinki\char'100science.unideb.hu; istvan.pink\char'100sbg.ac.at}

\address{V. Ziegler \newline
         \indent University of Salzburg \newline
         \indent Hellbrunnerstrasse 34/I \newline
         \indent A-5020 Salzburg, Austria}

\email{volker.ziegler\char'100sbg.ac.at}

\begin{document}

\begin{abstract}
In this paper, we find all integers $c$ having at least two representations as a difference between a Fibonacci number and a Tribonacci number.
\end{abstract}

\maketitle

\section{Introduction}

Pillai's famous conjecture first formulated in \cite{Pillai:1936} states that the Diophantine equation
\begin{equation}\label{eq:Pillai} a^x-b^y=c\end{equation}
has for any fixed integer $c>0$ at most finitely many solutions $a,b,x,y$ in positive integers.
This conjecture is still open for all $c\neq 1$. Note that
the case $c=1$ is Catalan's conjecture which has been solved by Mih{\v{a}}ilescu \cite{Mihailescu:2004}.
If we leave $a,b$ and $c$ fixed, then much more is known about
the solutions $(x,y)$. For instance Pillai~\cite{Pillai:1936} showed that if $c$ is larger
than some constant depending on $a$ and $b$, then Diophantine
equation \eqref{eq:Pillai}
has at most one solution. In particular, he conjectured that in the case that $a=3$ and $b=2$ Diophantine equation \eqref{eq:Pillai}
has at most one solution if $c>13$. This conjecture
was confirmed by Stroeker and Tijdeman \cite{Stroeker:1982} and their result was further improved by Bennett \cite{Bennett:2001},
who showed that for fixed $a,b$ and
$c$ equation \eqref{eq:Pillai} has at most two solutions.

Recently Ddamulira, Luca and Rakotomalala \cite{Luca15} considered the Diophantine equation
\begin{equation}\label{eq:Luca} F_n - 2^m = c, \end{equation}
where $c$ is a fixed integer and $\{F_n\}_{n \geqslant 0}$ is the sequence of Fibonacci numbers given
by $F_0 = 0$, $F_1 = 1$ and $F_{n+2} = F_{n+1} + F_n$ for all $n \geqslant 0$.
This type of equation can be seen as a variation of Pillai's equation. However Ddamulira et.al. proved
that the only integers $c$ having at least two representations of the form
$F_n - 2^m$ are contained in the set  $\mathcal C = \{ 0, 1, -3, 5, -11, -30, 85 \}$. Moreover, they
computed for all $c\in\mathcal C$ all representations of the from \eqref{eq:Luca}.

The purpose of this paper is to consider a related problem. Denote by $\{ T_m \}_{m\geqslant 0}$ the sequence of Tribonacci numbers given
by $T_0 = 0$, $T_1=T_2=1$ and $T_{m+3}=T_{m+2}+T_{m+1} + T_m$ for all $m \geqslant 0$. The main result of our paper is to find all nonzero
integers $c$ admitting at least two representations of the form $F_n - T_m$ for some positive integers $n$ and $m$. It is assumed that
representations with $n \in \lbrace 1, 2 \rbrace$ (for which $F_1=F_2=1$) as well as
representations with $m \in \lbrace 1, 2 \rbrace$ (for which $T_1=T_2$)
count as one representation to avoid trivial parametric families such as $1-1 = F_1-T_1 = F_2 - T_1=F_1-T_2=F_2-T_2$.
Therefore we assume that $n \geqslant 2$ and $m \geqslant 2$. We prove the following theorem:

\begin{theorem}\label{Main}
The only integers $c$ having at least two representations of the form $F_n - T_m$ come from the set
$$\mathcal C=\{ 0, 1, -1, -2,  -3, 4, -5, 6, 8, -10, 11, -11, -22, -23, -41, -60, -271 \}.$$
Furthermore, for each $c\in\mathcal C$ all representations of the form $c=F_n - T_m$ with integers $n \geqslant 2$ and $m \geqslant 2$ are:
\begin{equation*}
\begin{split}
0 & = 1 - 1 = 2 - 2 = 13 - 13 \; (= F_{2} - T_{2} = F_{3} - T_{3} = F_{7} - T_{6}),\\
1 & = 2 - 1 = 3 - 2 = 5 - 4 = 8 - 7 \\
& \qquad \qquad (= F_{3} - T_{2} = F_{4} - T_{3}= F_{5} - T_{4} = F_{6} - T_{5}),\\
-1 & = 1 - 2 = 3 - 4 \; (= F_{2} - T_{3} = F_{4} - T_{4}),\\
-2 & = 2 - 4 = 5 - 7 \; (= F_{3} - T_{4} = F_{5} - T_{5}),\\
-3 & = 1 - 4 = 21 - 24 \; (= F_{2} - T_{4} = F_{8} - T_{7}),\\
4 & = 5 - 1 = 8 - 4 \; (= F_{5} - T_{2} = F_{6} - T_{4}),\\
-5 & = 2 - 7 = 8 - 13 = 144 - 149 \; (= F_{3} - T_{5} = F_{6} - T_{6} = F_{12} - T_{10}),\\
6 & = 8 - 2 = 13 - 7 \; (= F_{6} - T_{3} = F_{7} - T_{5}),\\
8 & = 21 - 13 = 89 - 81 \; (= F_{8} - T_{6} = F_{11} - T_{9}),\\
-10 & = 3 - 13 = 34 - 44 \; (= F_{4} - T_{6} = F_{9} - T_{8}),\\
11 & = 13 - 2 = 55 - 44 \; (= F_{7} - T_{3} = F_{10} - T_{8}),\\
-11 & = 2 - 13 = 13 - 24 \; (= F_{3} - T_{6} = F_{7} - T_{7}),\\
-22 & = 2 - 24 = 121393 - 121415 \; (= F_{3} - T_{7} = F_{26} - T_{21}),\\
-23 & = 1 - 24 = 21 - 44 \; (= F_{2} - T_{7} = F_{8} - T_{8}),\\
-41 & = 3 - 44 = 233 - 274 \; (= F_{4} - T_{8} = F_{13} - T_{11}),\\
-60 & = 21 - 81 = 89 - 149 \; (= F_{8} - T_{9} = F_{11} - T_{10}),\\
-271 & = 3 - 274 = 233 - 504 \; (= F_{4} - T_{11} = F_{13} - T_{12}).
\end{split}
\end{equation*}
\end{theorem}


\section{Preliminaries}

In this section, the result of linear forms in logarithms by Baker and W\"ustholz~\cite{bawu93} is stated.
Besides, we state a lemma used by Ddamulira et.al. \cite{Luca15}, which is a slight variation of a result due to Dujella and Peth\H{o} \cite{dujella98},
of which is a generalization of a result due to Baker and Davenport \cite{BD69}. Both results will be used in the proof of Theorem \ref{Main}.

\subsection{A lower bound for linear forms in logarithms of algebraic numbers}

In 1993, Baker and W\"ustholz \cite{bawu93} obtained an explicit bound for linear forms in logarithms with a linear dependence on $\log B$.
It is a vast improvement compared with lower bounds with a dependence on higher powers of $\log B$ in preceding
publications by other mathematicians in particular Baker's original results \cite{Baker:1966}. The final structure
for the lower bound for linear forms in logarithms without an
explicit determination of the constant $C(k,d)$
involved has been established by W\"ustholz \cite{wu88} and the precise determination of that constant
is the central aspect of \cite{bawu93} (see also \cite{bawu07}). The improvement was mainly due to the use of the
analytic subgroup theorem established by W\"ustholz \cite{wu89}. We shall now state the result of Baker and W\"ustholz.

Denote by $\alpha_1, \dots, \alpha_k$ algebraic numbers, not $0$ or $1$, and by $\log \alpha_1, \dots, \log \alpha_k$
a fixed determination of their logarithms. Let $K=\bbQ(\alpha_1, \dotso, \alpha_k)$ and let $d=[K:\bbQ]$ be the degree of $K$ over $\bbQ$.
For any $\alpha \in K$, suppose that its minimal polynomial over the integers is
\[ g(x) = a_0 x^{\delta} + a_1 x^{\delta - 1} + \dotso + a_{\delta} = a_0 \prod_{j=1}^{\delta} (x - \alpha^{(j)}).\]
The absolute logarithmic Weil height of $\alpha$ is defined as
\[ h_0(\alpha) = \dfrac{1}{\delta} \left( \log |a_0| + \sum_{j=1}^{\delta} \log \left( \max \lbrace|\alpha^{(j)}|, 1 \rbrace \right)\right). \]
Then the modified height $h'(\alpha)$ is defined by
\[ h'(\alpha)=\dfrac{1}{d}\max \{h(\alpha), |\log\alpha|, 1\},\]
where $h(\alpha) = d h_0(\alpha)$ is the standard logarithmic Weil height of $\alpha$.

Let us consider the linear form
\[ L(z_1, \dotso, z_k)=b_1 z_1 + \dotso + b_k z_k, \]
where $b_1, \dotso, b_k$ are rational integers, not all 0 and define
\[ h'(L) = \dfrac{1}{d}\max\{h(L),1\}, \]
where $h(L) = d \log \left(\max_{1 \le j \le k} \left\{\dfrac{|b_j|}{b}\right\}\right)$ is the logarithmic Weil height of $L$,
where $b$ is the greatest common divisor of $b_1, \dotso, b_k$.
If we write $B=\max\lbrace|b_1|, \dotso, |b_k|, e \rbrace$, then we get
\[ h'(L)\leqslant \log B. \]

With these notations we are able to state the following result due to Baker and W\"ust\-holz~\cite{bawu93}.

\begin{theorem}\label{BaWu}
If ${\it\Lambda}=L(\log \alpha_1, \dotso, \log\alpha_k) \neq 0$, then
\[\log|{\it\Lambda}|\geqslant -C(k,d)h'(\alpha_1)\dotso h'(\alpha_k)h'(L), \]
where
\[ C(k,d)=18(k+1)!\,k^{k+1}(32d)^{k+2}\log(2kd). \]
\end{theorem}

With $|{\it\Lambda}| \leqslant \dfrac{1}{2}$, we have $\dfrac{1}{2}|{\it\Lambda}| \leqslant |{\it\Phi}| \leqslant 2|\it{\Lambda}|$, where
\[ {\it\Phi} = e^{{\it\Lambda}}-1 = \alpha_1^{b_1} \cdots \alpha_k^{b_k}-1, \]
so that
\[ \log|\alpha_1^{b_1} \cdots \alpha_k^{b_k}-1| \geqslant \log|{\it\Lambda}| - \log2. \]

We apply Theorem \ref{BaWu} mainly in a situation where $k=3$ and $d=6$. In this case we obtain
\[C(3,6)=18(4!) 3^4 (32\times 6)^5 (\log36)\approx 3.2718\ldots \times 10^{16}.\]
We will use this value throughout the paper without any further reference.


\subsection{A generalized result by Dujella and Peth\H{o}}
The following result will be used to reduce the huge upper bounds for $n$ and $m$ which appear during the course of the proof of
Theorem \ref{Main} (cf. Proposition \ref{prop:bound}). The following Lemma is stated in~\cite{Luca15}, which is regarded as a slight variation
of a result due to Dujella and Peth\H{o}~\cite{dujella98}, of which is a generalization of a result due to Baker and Davenport \cite{BD69}.
For a real number $x$, let $\|x\| = \min \{ |x-n| : n \in \bbZ \}$ be the distance from $x$ to the nearest integer.\

\begin{lemma} \label{Dujella}
Let $M$ be a positive integer, let $p/q$ be a convergent of the continued fraction of the irrational $\tau$ such that $q > 6M$,
and let $A, B, \mu$ be some real numbers with $A > 0$ and $B > 1$. Let $\varepsilon := \|\mu q\| - M \| \tau q\|$.
If $\varepsilon > 0$, then there is no solution to the inequality
\[ 0 < m \tau - n + \mu < AB^{-k}, \]
in positive integers $m, n$ and $k$ with
\[ m \leqslant M \quad \quad and \quad \quad k \geqslant \dfrac{\log (Aq/ \varepsilon )}{\log B}. \]
\end{lemma}


\section{Proof of Theorem \ref{Main}}
\subsection{Set up}

Assume that $(n,m) \neq (n_1,m_1)$ are pairs of indices such that
\begin{equation}\label{eqt}
F_n-F_{n_1} = T_m - T_{m_1}.
\end{equation}
We may assume that $m \neq m_1$, since otherwise $(n,m) = (n_1,m_1)$. Furthermore we assume that $m > m_1$.
Due to equation \eqref{eqt} and since the right hand side of equation \eqref{eqt} is positive,
we get that the left hand side of equation \eqref{eqt} is also positive and thus $n > n_1$. Therefore,
we have $n \geqslant 3$, $n_1 \geqslant 2$ and $m \geqslant 3$, $m_1 \geqslant 2$.

During the proof of Theorem \ref{Main} we use the Binet formulae for the Fibonacci sequence and Tribonacci sequence in the following form:

\textbf{Fibonacci sequence:}
\[ F_k = \dfrac{\alpha^k-\beta^k}{\alpha - \beta} \quad \mbox{ for all } k \geqslant 0, \]
where $\alpha=\frac{1+\sqrt{5}}{2}$ and $\beta =\frac{1-\sqrt{5}}{2}$ are the roots of the characteristic equation $x^2-x-1=0$. Besides, the inequality
\[ \alpha^{k-2} \leqslant F_k \leqslant \alpha^{k-1}\]
holds for all $k \geqslant 1$.

\textbf{Tribonacci sequence:}
\[ T_k = c_{\alpha} \alpha_T^k + c_{\beta} \beta_T^k + c_{\gamma} \gamma_T^k \quad \mbox{ for all } k \geqslant 0, \]
where $\alpha_T$, $\beta_T$ and  $\gamma_T$ are the roots of the characteristic equation $x^3-x^2 - x -1=0$, with
\begin{align*}
\alpha_T = & \ \dfrac{1}{3}\left( 1+\sqrt[3]{19+3\sqrt{33}} + \sqrt[3]{19-3\sqrt{33}} \right),\\
\beta_T = & \ \dfrac{1}{6}\left( 2- \sqrt[3]{19+3\sqrt{33}} - \sqrt[3]{19-3\sqrt{33}} \right)+\dfrac{\sqrt{3}}{6}i\left(\sqrt[3]{19+3\sqrt{33}}- \sqrt[3]{19-3\sqrt{33}} \right),\\
\gamma_T = & \ \dfrac{1}{6}\left( 2- \sqrt[3]{19+3\sqrt{33}} - \sqrt[3]{19-3\sqrt{33}} \right)-\dfrac{\sqrt{3}}{6}i\left(\sqrt[3]{19+3\sqrt{33}}- \sqrt[3]{19-3\sqrt{33}} \right),
\end{align*}
and the coefficients
\begin{align*}
 c_\alpha = &\ \dfrac{\alpha_T}{(\alpha_T-\beta_T)(\alpha_T-\gamma_T)}=\dfrac{1}{-\alpha_T^2+4\alpha_T-1},\\
 c_\beta= &\ \dfrac{\beta_T}{(\beta_T-\alpha_T)(\beta_T-\gamma_T)}=\dfrac{1}{-\beta_T^2+4\beta_T-1},\\
 c_\gamma =&\ \dfrac{\gamma_T}{(\gamma_T-\alpha_T)(\gamma_T-\beta_T)}= \dfrac{1}{-\gamma_T^2+4\gamma_T-1}
\end{align*}
are the roots of the polynomial $44x^3-2x -1$. Note that
\begin{align*}
& 1.839 < \alpha_T < 1.840 && 0.336 < c_{\alpha} < 0.337\\
& \beta_T= \overline{\gamma_T} && 0.737< |\beta_T|= |\gamma_T| < 0.738 \\
& c_{\beta} = \overline{c_{\gamma}}  && 0.259< |c_{\beta}| = |c_{\gamma}| < 0.260.
\end{align*}

Finally let us state several useful inequalities. For instance
\[ \alpha_T^{k-2} \leqslant T_k \leqslant \alpha_T^{k-1} \quad \text{for all } k \geqslant 1. \]
which was already shown in \cite{Luca12}. Using equation \eqref{eqt} we get that
\begin{equation}\label{ineq1}
\alpha^{n-4} \leqslant F_{n-2} \leqslant F_n-F_{n_1} = T_m - T_{m_1} < T_m \leqslant \alpha_T^{m-1},
\end{equation}
and similarly we get
\begin{equation}\label{ineq2}
\begin{split}
\alpha^{n-1} \geqslant F_n > F_n-F_{n_1} &= T_m - T_{m_1} \geqslant T_m - T_{m-1} = \\
& T_{m-2}+T_{m-3} \geqslant \alpha_T^{m-4}+\alpha_T^{m-5} > 2.83 \alpha_T^{m-5}.
\end{split}
\end{equation}
Thus
\begin{equation}\label{computer}
n-4< \dfrac{\log \alpha_T}{\log \alpha} (m-1) \quad \mbox{and} \quad n-3 > \dfrac{\log \alpha_T}{\log \alpha} (m-5),
\end{equation}
where $\dfrac{\log \alpha_T}{\log \alpha} \approx 1.2663\dots$.

Inequality \eqref{computer} implies that if $n < 300$, then $m < 240$. By a brute force computer
enumeration for $2 \leqslant n_1 < n < 300$ and $2 \leqslant m_1 < m < 240$
we found all solutions listed in Theorem \ref{Main}. Thus we may assume from now on that $n\geq 300$.


\subsection{Linear forms in logarithms}

Since $n \ge 300$, by the first inequality of \eqref{computer} we obtain that $m \ge 235$
which combined with the second inequality of \eqref{computer} implies that $n>m$.
Moreover, we have
$$ \dfrac{\alpha^n - \beta^n}{\sqrt{5}} -\dfrac{\alpha^{n_1} - \beta^{n_1}}{\sqrt{5}} =
\left( c_{\alpha} \alpha_T^m + c_{\beta} \beta_T^m + c_{\gamma} \gamma_T^m \right) -
\left( c_{\alpha} \alpha_T^{m_1} + c_{\beta} \beta_T^{m_1} + c_{\gamma} \gamma_T^{m_1} \right).$$
Collecting the ``large'' terms on the left hand side of the equation we obtain
\begin{align*}
\left| \dfrac{\alpha^n}{\sqrt{5}} - c_{\alpha} \alpha_T^m \right| = &\ \left| \dfrac{\beta^n}{\sqrt{5}}
+ \dfrac{\alpha^{n_1} - \beta^{n_1}}{\sqrt{5}} + \left(  c_{\beta} \beta_T^m + c_{\gamma} \gamma_T^m \right)
- \left( c_{\alpha} \alpha_T^{m_1} + c_{\beta} \beta_T^{m_1} + c_{\gamma} \gamma_T^{m_1} \right) \right|\\
\leqslant &\ \dfrac{\alpha^{n_1}}{\sqrt{5}}  + c_{\alpha}\alpha_T^{m_1} + \dfrac{\left| \beta\right|^n}{\sqrt{5}}
+ \dfrac{\left| \beta\right|^{n_1}}{\sqrt{5}} +  \left|  c_{\beta}\right| \left| \beta_T \right|^m +
\left| c_{\gamma} \right| \left| \gamma_T \right|^m \\
&\ +\left| c_{\beta}\right| \left| \beta_T\right|^{m_1} + \left| c_{\gamma}\right| \left|  \gamma_T \right|^{m_1} \\
< &\ \dfrac{\alpha^{n_1}}{\sqrt{5}} + c_{\alpha} \alpha_T^{m_1} + 0.46 \\
< &\ 0.92 \max \{ \alpha^{n_1},\alpha_T^{m_1} \}.
\end{align*}
Dividing by $c_{\alpha} \alpha_T^m$ we get
\begin{align*}
\left| (\sqrt{5}c_{\alpha})^{-1} \alpha^n \alpha_T^{-m} - 1 \right| &<
\max \left\{ \dfrac{0.92}{c_\alpha \alpha_T^{m} } \alpha^{n_1},\dfrac{0.92}{c_\alpha} \alpha_T^{m_1-m} \right\}\\
 &<  \max \left\{ 2.74 \dfrac{\alpha^{n_1}}{\alpha_T} \dfrac{1}{\alpha^{n-4}},2.74 \alpha_T^{m_1-m} \right\}.
\end{align*}
Hence we obtain the inequality
\begin{equation}\label{Case0}
\left| (\sqrt{5}c_{\alpha})^{-1} \alpha^n \alpha_T^{-m} - 1 \right|< \max \{ \alpha^{n_1-n+5},\alpha_T^{m_1-m+2} \}.
\end{equation}

Let us introduce
$${\it\Lambda} = n\log \alpha - m \log \alpha_T - \log (\sqrt{5}c_{\alpha})$$
and assume that $|{\it\Lambda}| \leqslant 0.5$.
Further, we put
$${\it\Phi} = e^{\it\Lambda}-1 = (\sqrt{5}c_{\alpha})^{-1} \alpha^n \alpha_T^{-m} - 1$$
and use the theorem of Baker and W\"ustholz (Theorem \ref{BaWu}) with the data
\[k = 3, \quad \alpha_1 = \sqrt{5} c_\alpha, \quad b_1 = -1, \quad \alpha_2 = \alpha, \quad b_2 = n, \quad \alpha_3 = \alpha_T,  \quad b_3 = -m. \]
With this data we have $K = \bbQ(\sqrt{5}, \alpha_T)$, i.e. $d = 6$, and $B=n$. Notice that the minimal polynomial of $\alpha_1$ is
$1936 x^6 - 880 x^4 +100 x^2-125$, and we conclude that $h'(\alpha_1) = \frac{1}{6}\log1936$. Further we obtain by a simple computation that
$h'(\alpha_2) = \frac{1}{2} \log \alpha$ and $h'(\alpha_3) = \frac{1}{3}\log \alpha_T$.

Before we can apply Theorem \ref{BaWu} we have to show that ${\it\Phi} \neq 0$. Assume to the contrary that ${\it\Phi} = 0$,
then $\alpha^{2n} = 5 c_\alpha^2 \alpha_T^{2m}$. But $\alpha^{2n} \in \bbQ(\sqrt{5})\setminus \bbQ$ whereas
$5 c_\alpha^2 \alpha_T^{2m} \in \bbQ(\alpha_T)$. Thus ${\it\Phi} = 0$ is impossible due to the fact that $\bbQ(\sqrt{5})\cap\bbQ(\alpha_T)=\bbQ$.

Apply Theorem \ref{BaWu} yields
\begin{equation*}
\log |{\it\Phi}| \geqslant -C(3,6) \left( \dfrac{1}{6}\log1936 \right)
\left( \dfrac{1}{2}  \log \alpha \right) \left( \dfrac{1}{3} \log \alpha_T \right) \log n - \log 2
\end{equation*}
and together with inequality \eqref{Case0} we have
\[\min \{ (n-n_1-5)\log \alpha,(m-m_1-2)\log\alpha_T \} \leqslant 2.02 \times 10^{15} \log n. \]
Thus we have proved so far:

\begin{lemma}\label{lem:Case0}
 Assume that $(n,m,n_1,m_1)$ is a solution to equation \eqref{eqt} with $m>m_1$. Then we have
 $$\min \{ (n-n_1)\log \alpha,(m-m_1)\log\alpha_T \} < 2.03 \times 10^{15} \log n.$$
\end{lemma}

Note that in the case that $|{\it\Lambda}| > 0.5$ inequality \eqref{Case0} is possible only if
either $n-n_1 \leqslant 5$ or $m-m_1 \leqslant 2$, which is covered by the bound provided by Lemma \ref{lem:Case0}.

Now we have to distinguish between the following two cases:

\noindent\textbf{Case 1.} Let us assume that
$$\min\lbrace (n-n_1)\log \alpha, (m-m_1)\log \alpha_T \rbrace = (n-n_1) \log \alpha.$$
We rewrite equation \eqref{eqt} as
\begin{align*}
\left| \dfrac{\alpha^n-\alpha^{n_1}}{\sqrt{5}} - c_{\alpha} \alpha_T^m \right| = &\
\left|  -c_{\alpha} \alpha_T^{m_1} + \dfrac{\beta^n}{\sqrt{5}} - \dfrac{\beta^{n_1}}{\sqrt{5}} + \right.\\
& \quad \left. \phantom{\dfrac{\beta^{n_1}}{\sqrt{5}}} \left(  c_{\beta} \beta_T^m + c_{\gamma} \gamma_T^m \right) - \left( c_{\beta} \beta_T^{m_1} +  c_{\gamma} \gamma_T^{m_1} \right) \right|\\
 \leqslant &\  c_{\alpha} \alpha_T^{m_1} + \dfrac{\left| \beta^n \right| }{\sqrt{5}} + \dfrac{\left|  \beta^{n_1} \right|}{\sqrt{5}}+
 \left| c_{\beta} \right| \left| \beta_T^m \right| + \left|  c_{\gamma} \right| \left| \gamma_T^m \right| + \\
 &\ \left| c_{\beta}\right| \left| \beta_T^{m_1} \right| +  \left| c_{\gamma}\right| \left|  \gamma_T^{m_1} \right|\\
\end{align*}
and obtain that
$$\left| \dfrac{\alpha^{n-n_1}-1}{\sqrt{5}} \alpha^{n_1} - c_{\alpha} \alpha_T^m \right|  <  (c_{\alpha} + 0.14) \alpha_T^{m_1}.$$
Dividing by $c_{\alpha} \alpha_T^{m}$ we get the inequality
\begin{equation} \label{Case1}
\left| \dfrac{\alpha^{n-n_1}-1}{\sqrt{5}c_{\alpha}} \alpha^{n_1} \alpha_T^{-m} - 1 \right| < 1.42\alpha_T^{m_1-m}.
\end{equation}

\noindent\textbf{Case 2.} Let us assume that
$$\min\lbrace (n-n_1)\log \alpha, (m-m_1)\log \alpha_T \rbrace = (m-m_1)\log \alpha_T.$$
We rewrite equation \eqref{eqt} as
\begin{align*}
\left| \dfrac{\alpha^n }{\sqrt{5}} -c_{\alpha} \alpha_T^{m} + c_{\alpha} \alpha_T^{m_1} \right|  = & \
\left|\dfrac{\beta^n}{\sqrt{5}}+\dfrac{\alpha^{n_1} - \beta^{n_1}}{\sqrt{5}} + c_{\beta} \beta_T^m + c_{\gamma} \gamma_T^m - c_{\beta} \beta_T^{m_1} - c_{\gamma} \gamma_T^{m_1} \right|\\
 \leqslant &\ \dfrac{\left|\beta^n \right|}{\sqrt{5}} +\dfrac{\alpha^{n_1}}{\sqrt{5}} + \dfrac{ \left|\beta^{n_1} \right|}{\sqrt{5}}
 + \left| c_{\beta}\right| \left| \beta_T^m \right| + \left| c_{\gamma} \right| \left|  \gamma_T^m \right|\\
 &\ + \left|  c_{\beta} \right| \left|  \beta_T^{m_1} \right| + \left|  c_{\gamma} \right| \left| \gamma_T^{m_1} \right|.
\end{align*}
Thus we get
$$\left|\alpha^n -\sqrt{5} c_{\alpha}(\alpha_T^{m-m_1} -1)\alpha_T^{m_1}  \right|  < 1.4 \alpha^{n_1}.$$
Dividing both sides by $\sqrt{5} c_{\alpha}(\alpha_T^{m-m_1} -1)\alpha_T^{m_1}$ we get by using inequality \eqref{ineq1} the following inequality:
\begin{equation} \label{Case2}
\left|\frac{\alpha^n \alpha_T^{-m_1}}{\sqrt{5}\, c_{\alpha}(\alpha_T^{m-m_1} -1)}  -1 \right| <
\dfrac{1.4 }{\sqrt{5}\,c_{\alpha} (1 -\alpha_T^{m_1-m})\alpha_T} \dfrac{\alpha^{n_1}}{\alpha_T^{m-1}} < 2.22 \alpha^{n_1-n+4}.
\end{equation}

We want to apply Theorem \ref{BaWu} to both inequalities \eqref{Case1} and \eqref{Case2} respectively.
Let us consider the first case more closely. We write
$${\it\Lambda}_1 =n_1 \log \alpha - m \log \alpha_T + \log \left( \dfrac{\alpha^{n-n_1}-1}{\sqrt{5}c_{\alpha}} \right)$$
and assume that
$|{\it\Lambda}_1| \leqslant 0.5$. Further, we put
$${\it\Phi}_1 = e^{{\it\Lambda}_1}-1= \dfrac{\alpha^{n-n_1}-1}{\sqrt{5}c_{\alpha}} \alpha^{n_1} \alpha_T^{-m} - 1$$
and aim to apply Theorem \ref{BaWu} by taking $K = \bbQ(\sqrt{5}, \alpha_T)$, i.e. $d = 6$, $k=3$ and $B = n$. Further, we have
\[ \alpha_1= \dfrac{\alpha^{n-n_1}-1}{\sqrt{5}c_{\alpha}}, \quad b_1 = 1, \quad \alpha_2 =  \alpha, \quad b_2 = n_1 ,
\quad \alpha_3 = \alpha_T, \quad b_3=-m.\]
Let us estimate the height of $\alpha_1$. Notice that $h(\alpha_1) \leqslant h(\eta_1)+h(\eta_2)$,
where $\eta_1 = \frac{\alpha^{n-n_1} - 1}{\sqrt{5}} $ and $\eta_2 = \frac{1}{c_\alpha }$.
The minimal polynomial of $\eta_1$ divides (e.g. see \cite{Luca15})
\[ 5X^2-5F_{n-n_1}X - ((-1)^{n-n_1}+1-L_{n-n_1}),\]
where $\lbrace L_k \rbrace_{k\geqslant 0}$ is the Lucas companion sequence of the Fibonacci sequence given by $L_0 = 2, L_1 = 1, L_{k+2}=L_{k+1}+L_k$
for $k \geqslant 0$. Its Binet formula for the general term is $L_k = \alpha^k + \beta^k$ for all $k\geqslant 0$. Thus (cf. \cite{Luca15}),
\[ h_0(\eta_1)\leqslant \dfrac{1}{2} \left( \log 5 + \log\left( \dfrac{\alpha^{n-n_1}+1}{\sqrt{5}} \right) \right). \]
Thus Lemma \ref{lem:Case0} yields an upper bound
$$ h_0(\eta_1) < \dfrac{1}{2} \log\left(2\sqrt{5} \alpha^{n-n_1} \right) < \dfrac{1}{2} (n-n_1+4) \log \alpha < 1.02 \times 10^{15}\log n,$$
i.e. $h(\eta_1) < 6\times 1.02 \times 10^{15}\log n$. Since $h_0(\eta_2) = h_0(c_\alpha) = \frac{1}{3}\log 44$, i.e. $h(\eta_2) = 2\log44$,
we have $h(\alpha_1) \leqslant 6\times 1.02 \times 10^{15}\log n + 2\log44$ and finally we obtain that
\[ h'(\alpha_1) < 1.03 \times 10^{15}\log n. \]
Moreover, we have that $h'(\alpha_2) = \frac{1}{2}\log \alpha$ and $h'(\alpha_3) = \frac{1}{3}\log \alpha_T$ as before.

Now let us turn to the second case. We write
$${\it\Lambda}_2 =n \log \alpha - m_1 \log \alpha_T - \log \left( \sqrt{5}\,c_{\alpha} (\alpha_T^{m-m_1}-1) \right)$$
and assume that $|{\it\Lambda}_2| \leqslant 0.5$. Further, we put
$${\it\Phi}_2 = e^{{\it\Lambda}_2}-1 = (\sqrt{5}\,c_{\alpha}(\alpha_T^{m-m_1} -1))^{-1}\alpha^n \alpha_T^{-m_1} -1$$
and aim to apply Theorem \ref{BaWu}. As in the previous case we take $K = \bbQ(\sqrt{5}, \alpha_T)$, i.e. $d = 6$, $k=3$ and $B = n$. Further, we have
\[ \alpha_1= \sqrt{5}\, c_{\alpha}(\alpha_T^{m-m_1} -1), \quad b_1 = -1, \quad \alpha_2 = \alpha, \quad b_2 = n , \quad \alpha_3 = \alpha_T, \quad b_3=-m_1.\]
Again, we have to estimate $h(\alpha_1)$ and therefore note that $h(\alpha_1) \leqslant h(\eta_1)+h(\eta_2)+h(\eta_3)$, where $\eta_1 = \alpha_T^{m-m_1} -1$, $\eta_2 = c_{\alpha}$ and $\eta_3 = \sqrt{5}$.
By applying Lemma \ref{lem:Case0} we get
\begin{align*}
h_0(\eta_1) & \leqslant h_0(\alpha_T^{m-m_1}) +h_0(-1) + \log2\\
& = (m-m_1)h_0(\alpha_T) + \log2 = \dfrac{m-m_1}{3}\log \alpha_T + \log 2\\
& < \dfrac{1}{3} \times 2.03 \times 10^{15} \log n + \log 2 .
\end{align*}
Thus
$$h(\alpha_1) < 6 \left( \dfrac{1}{3} \times 2.03 \times 10^{15} \log n + \log 2 + \dfrac{1}{3} \log 44 + \log \sqrt{5} \right)$$
and therefore
\begin{align*}
h'(\alpha_1) < 6.77 \times 10^{14} \log n<1.03 \times 10^{15} \log n.
\end{align*}
Once again, we have that $h'(\alpha_2) = \frac{1}{2}\log \alpha$ and $h'(\alpha_3) = \frac{1}{3}\log \alpha_T$.

In particular, we have shown in both cases that
$$
h'(\alpha_1) < 1.03 \times 10^{15} \log n, \quad h'(\alpha_2) = \dfrac{1}{2}\log \alpha, \quad h'(\alpha_3) = \dfrac{1}{3}\log \alpha_T, \quad B = n.
$$

But, before we can apply Theorem \ref{BaWu} we have to ensure that ${\it\Phi}_i \neq 0$ for $i=1, 2$.
Firstly we deal with the assumption that ${\it\Phi}_1 = 0$,
i.e. $\alpha^{n}-\alpha^{n_1} = \sqrt{5}c_\alpha \alpha_T^m $. This is impossible if
$\sqrt{5}c_\alpha \alpha_T^m \in \bbQ(\sqrt{5}, \alpha_T)$ but $\notin \bbQ(\sqrt{5})$.
Therefore let us assume that $\sqrt{5}c_\alpha \alpha_T^m \in \bbQ(\sqrt{5})$, hence
$\sqrt{5}c_\alpha \alpha_T^m=y \sqrt{5}$ for some $y \in \bbQ$. Let
$\sigma\neq \mathrm{id}$ be the unique non-trivial $\bbQ$-automorphism over $\bbQ(\sqrt 5)$. Then we get
\[\alpha^{n}-\alpha^{n_1} =\sqrt{5}c_\alpha \alpha_T^m=y\sqrt{5}=-\sigma(\sqrt{5}c_\alpha \alpha_T^m) =
- \sigma(\alpha^{n}-\alpha^{n_1})= \beta^{n_1}-\beta^{n}.\]
However, the absolute value of $\alpha^{n}-\alpha^{n_1}$ is at least
$\alpha^{n}-\alpha^{n_1} \geqslant \alpha^{n-2} \geqslant \alpha^{298}>2$
whereas the absolute value of $\beta^{n_1}-\beta^n $ is at most
$|\beta^{n_1}-\beta^n| \leqslant |\beta|^{n_1} + |\beta|^n < 2$. By this obvious contradiction we conclude that ${\it\Phi}_1 \neq 0$.

Now let us consider the second case and assume for the moment that ${\it\Phi}_2 = 0$, i.e. $\alpha^{2n} = 5 \alpha_T^{2m_1} c_\alpha^2 (\alpha_T^{m-m_1} - 1)^2$.
However, $\alpha^{2n} \in \bbQ(\sqrt{5}) \setminus \bbQ$, whereas $ 5 \alpha_T^{2m_1} c_\alpha^2 (\alpha_T^{m-m_1} - 1)^2 \in \bbQ(\alpha_T)$.
Thus we obtain also in this case a contradiction and we also conclude in this case that ${\it\Phi}_2 \neq 0$.

Now, we are ready to apply Theorem \ref{BaWu} and get
\begin{align*}
\log |{\it\Phi}_i|  > & -C(3,6) \left( 1.03 \times 10^{15} \log n \right)
 \left(  \dfrac{1}{2}\log \alpha \right) \left( \dfrac{1}{3}\log \alpha_T \right) \log n - \log 2\\
> & -1.65 \times 10^{30} (\log n)^2
\end{align*}
for $i = 1,2$. Combining this inequality with the inequalities \eqref{Case1} and \eqref{Case2}, we obtain
\[ (m_1-m)\log \alpha_T + \log 1.42 > -1.65 \times 10^{30} (\log n)^2 \]
and
\[(n_1-n+4)\log \alpha + \log 2.22 > -1.65 \times 10^{30} (\log n)^2 \]
respectively. These two inequalities yield together with Lemma \ref{lem:Case0} the following lemma:

\begin{lemma}\label{lem:Case1-2}
 Assume that $(n,m,n_1,m_1)$ is a solution to equation \eqref{eqt} with $m>m_1$. Then we have
 $$\max \{ (n-n_1)\log \alpha,(m-m_1)\log\alpha_T \} < 1.66 \times 10^{30} (\log n)^2.$$
\end{lemma}

Note that in view of inequality \eqref{Case1} $|{\it\Lambda}_1| > 0.5$ is possible only if $m-m_1 =1$
and in view of inequality \eqref{Case2} $|{\it\Lambda}_2| > 0.5$
is possible only if $n-n_1 \leqslant 6$ respectively. Both cases are covered by the bound provided by Lemma \ref{lem:Case1-2}.

One more time we have to apply Theorem \ref{BaWu}. This time we rewrite equation \eqref{eqt} as
\begin{multline*}
\left| \dfrac{\alpha^n}{\sqrt{5}}  -\dfrac{\alpha^{n_1}}{\sqrt{5}} - c_{\alpha} \alpha_T^m  + c_{\alpha} \alpha_T^{m_1} \right| \\
=\left| \dfrac{\beta^n}{\sqrt{5}} - \dfrac{\beta^{n_1}}{\sqrt{5}}+ c_{\beta} \beta_T^m + c_{\gamma} \gamma_T^m - c_{\beta} \beta_T^{m_1} - c_{\gamma} \gamma_T^{m_1} \right| < 0.46
\end{multline*}
Dividing both sides by $c_{\alpha} \alpha_T^{m_1}(\alpha_T^{m-m_1}-1)$ we get by applying inequality \eqref{ineq1}
\begin{equation} \label{Case3}
 \left| \dfrac{\alpha^{n-n_1}-1}{\sqrt{5}\,c_{\alpha}(\alpha_T^{m-m_1}-1)} \alpha^{n_1}  \alpha_T^{-m_1} - 1 \right| <
 \dfrac{0.46}{c_{\alpha}(1- \alpha_T^{m_1-m})\alpha_T} \dfrac{1}{\alpha_T^{m-1}} < 1.64 \alpha^{4-n}.
\end{equation}
In this final step we consider the linear form
$${\it\Lambda}_3 =n_1 \log \alpha - m_1 \log \alpha_T + \log \left( \dfrac{\alpha_T^{n-n_1}-1}{\sqrt{5}\,c_{\alpha} (\alpha_T^{m-m_1}-1)} \right)$$
and assume that $|{\it\Lambda}_3| \leqslant 0.5$. Further, we put
$${\it\Phi}_3 = e^{{\it\Lambda}_3} - 1 =\dfrac{\alpha^{n-n_1}-1}{\sqrt{5}\,c_{\alpha}(\alpha_T^{m-m_1}-1)} \alpha^{n_1}  \alpha_T^{-m_1} - 1.$$
As before we take $K = \bbQ(\sqrt{5}, \alpha_T)$, i.e. $d = 6$, $k=3$, $B=n$ and we have
\[ \alpha_1 = \dfrac{\alpha^{n-n_1}-1}{\sqrt{5}\,c_{\alpha}(\alpha_T^{m-m_1}-1)}, \quad b_1 = 1, \quad \alpha_2 = \alpha, \quad b_2 = n_1, \quad \alpha_3 = \alpha_T, \quad b_3 = -m_1.\]
By Lemma \ref{lem:Case1-2} and similar computations as done before we obtain that
\[
h\left(\dfrac{\alpha^{n-n_1}-1}{\sqrt{5}c_{\alpha}}\right) \leqslant
6 \left( \dfrac{1}{2}(n-n_1+4) \log \alpha \right) + 2\log 44 < 3 \times \left( 1.67 \times 10^{30} (\log n)^2\right) \]
and
\[ h(\alpha_T^{m-m_1}-1) \leqslant 6 \left( \dfrac{m-m_1}{3}\log \alpha_T + \log 2 \right) < 2 \times \left(1.67 \times 10^{30}(\log n)^2\right).\]
Therefore we find the upper bound
$$h(\alpha_1) \leqslant h(\dfrac{\alpha^{n-n_1}-1}{\sqrt{5}c_{\alpha}} ) + h(\alpha_T^{m-m_1}-1) < 5 \times \left(1.67 \times 10^{30}(\log n)^2\right)$$
and thus
\[ h'(\alpha_1) < \dfrac{5}{6} \times \left(1.67  \times 10^{30}(\log n)^2\right). \]
As before, we have $h'(\alpha_2) = \dfrac{1}{2}\log \alpha$ and $h'(\alpha_3) = \dfrac{1}{3}\log \alpha_T$.

Using similar arguments as in the proof that ${\it\Phi}_1 \neq 0$ we can show that ${\it\Phi}_3 \neq 0$. Now an application of Theorem \ref{BaWu} yields
$$
\log |{\it\Phi}_3|
 > -C(3,6)\left(\dfrac{5}{6} \times 1.67 \times 10^{30} (\log n)^2 \right)
   \left( \dfrac{1}{2}\log \alpha \right) \left( \dfrac{1}{3}\log \alpha_T \right) \log n - \log 2.
$$

Combining this inequality with inequality \eqref{Case3} we get
$$(n-4) \log \alpha < 2.23 \times 10^{45} (\log n)^3$$
which yields $n < 8 \times 10^{51}$.

Similarly as in the cases above the assumption that $|{\it\Lambda}_3| > 0.5$ leads in view of inequality \eqref{Case3} to $n \leqslant 5$.
Let us summarize the results of this subsection:

\begin{proposition}\label{prop:bound}
 Assume that $(n,m,n_1,m_1)$ is a solution to equation \eqref{eqt} with $m>m_1$. Then we have that $n <8 \times 10^{51}$.
\end{proposition}

\begin{remark}
The theorem of Baker and W\"ustholz (Theorem \ref{BaWu}) can be easily applied. However, a slightly sharper bound for $n$, namely $n<6\times 10^{48}$,
may be obtained if one uses Matveev's result \cite{matveev00} instead. However, the improvement is not crucial in view of our next step,
the reduction of our upper bound for $n$ using the method of Baker and Davenport (Lemma \ref{Dujella}).
\end{remark}


\subsection{Generalized method of Baker and Davenport}
In our final step we reduce the huge upper bound for $n$ form Proposition \ref{prop:bound} by applying several times Lemma \ref{Dujella}.
In this subsection we follow the ideas from \cite{Luca15}. First, we consider inequality \eqref{Case0} and recall that
\[ {\it\Lambda}= n\log \alpha - m \log \alpha_T - \log (\sqrt{5}c_\alpha). \]
For technical reasons we assume that $\min\{n-n_1, m-m_1\} \geqslant 20$. In the case that
this condition fails we consider one of the following inequalities instead:

\begin{itemize}
\item if $n-n_1 < 20$ but $m-m_1 \geqslant 20$, we consider inequality \eqref{Case1};
\item if $n-n_1 \geqslant 20$ but $m-m_1 < 20$, we consider inequality \eqref{Case2};
\item if both $n-n_1 < 20$ and $m-m_1 < 20$, we consider inequality \eqref{Case3}.
\end{itemize}

Let us start by considering inequality \eqref{Case0}. Since we assume that $\min\{n-n_1, m-m_1\} \geqslant 20$
we get $|{\it\Phi}| = |e^{\it\Lambda}-1|<\frac{1}{4}$, hence $|{\it\Lambda}|< \frac{1}{2}$.
And, since $|x| < 2|e^x-1|$ holds for all $x \in (-\frac{1}{2}, \frac{1}{2})$ we get
\[ |{\it\Lambda}|<2\max \{ \alpha^{n_1 - n+5}, \alpha_T^{m_1 - m + 2} \} \leqslant \max \{ \alpha^{n_1 - n+7}, \alpha_T^{m_1 - m + 4} \}. \]
Assume that ${\it\Lambda} > 0$. Then we have the inequality
\begin{equation*}
\begin{split}
0< n \left( \dfrac{\log \alpha}{\log \alpha_T} \right) -m + \dfrac{\log(1/(\sqrt{5}c_\alpha))}{\log \alpha_T} < &
\max\left\{ \dfrac{\alpha^7}{\log \alpha_T} \alpha^{-(n-n_1)}, \dfrac{\alpha_T^4}{\log \alpha_T} \alpha_T^{-(m-m_1)} \right\}\\
< & \max\left\{ 48 \alpha^{-(n-n_1)},19 \alpha_T^{-(m-m_1)}\right\}
\end{split}
\end{equation*}
and we apply Lemma \ref{Dujella} with
\[ \tau=\dfrac{\log \alpha}{\log \alpha_T}, \quad \mu = \dfrac{\log(1/(\sqrt{5}c_\alpha))}{\log \alpha_T}, \quad (A,B)=(48, \alpha) \mbox{ or } (19, \alpha_T). \]
Let $\tau = [a_0, a_1, a_2, \dotso]=[0,1,3,1,3,13,2,1,8,3,1,5,\dotso]$ be the continued fraction of $\tau$.
Moreover, we choose $M = 8 \times 10^{51}$ and consider the $104$-th convergent
\[\dfrac{p}{q} = \dfrac{p_{104}}{q_{104}}=
\dfrac{528419636478855291192208008138409657842309076397924033}{669159011284129920139468279297504453112608160771671810}, \]
with $q=q_{104} > 6M$. This yields $\varepsilon > 0.068$ and therefore either
\[ n-n_1 \leqslant \dfrac{\log(48q/0.068)}{\log \alpha} < 272, \quad \mbox{or} \quad m-m_1 \leqslant \dfrac{\log(19q/0.068)}{\log \alpha_T} < 213. \]
Thus, we have either $n-n_1 \leqslant 271$, or $m-m_1 \leqslant 212$.

In the case of ${\it\Lambda}<0$ we consider the following inequality:
\begin{equation*}
\begin{split}
 0< m \left( \dfrac{\log \alpha_T}{\log \alpha} \right) -n + \dfrac{\log(\sqrt{5}c_\alpha)}{\log \alpha} < &
\max\left\{ \dfrac{\alpha^7}{\log \alpha_T} \alpha^{-(n-n_1)}, \dfrac{\alpha_T^4}{\log \alpha_T} \alpha_T^{-(m-m_1)} \right\}\\
<& \max\left\{ 61 \alpha^{-(n-n_1)},24 \alpha_T^{-(m-m_1)}\right\}
\end{split}
\end{equation*}
instead and apply Lemma \ref{Dujella} with
\[ \tau=\dfrac{\log \alpha_T}{\log \alpha}, \quad \mu = \dfrac{\log(\sqrt{5}c_\alpha)}{\log \alpha}, \quad (A,B)=(61, \alpha) \mbox{ or } (24, \alpha_T). \]
Let $\tau = [a_0, a_1, a_2, \dotso]=[1, 3, 1, 3, 13, 2, 1, 8, 3, 1, 5, 2, \dotso]$ be the continued fraction of $\tau$.
Again, we choose $M = 8 \times 10^{51}$ but in this case we consider instead of the $104$-th convergent the $103$-rd convergent
\[\dfrac{p}{q} = \dfrac{p_{103}}{q_{103}}=
\dfrac{669159011284129920139468279297504453112608160771671810}{528419636478855291192208008138409657842309076397924033}, \]
with $q > 6M$. This yields $\varepsilon > 0.067$ and again we obtain either
\[ n-n_1 \leqslant \dfrac{\log(61q/0.067)}{\log \alpha} < 272, \quad \mbox{or} \quad m-m_1 \leqslant \dfrac{\log(24q/0.067)}{\log \alpha_T} < 213. \]
These bounds agree with the bounds obtained in the case that ${\it\Lambda} > 0$. As a conclusion,
we have either $n-n_1 \leqslant 271$ or $m-m_1 \leqslant 212$ whenever ${\it\Lambda} \neq 0$.

Now, we have to distinguish between the two cases $n-n_1 \leqslant 271$ and $m-m_1 \leqslant 212$.
First, let us assume that $n-n_1 \leqslant 271$. In this case we consider inequality~\eqref{Case1} and assume that $m-m_1 \geqslant 20$. Recall that
\[ {\it\Lambda}_1= n_1\log \alpha - m \log \alpha_T + \log\left( \dfrac{\alpha^{n-n_1}-1}{\sqrt{5}c_\alpha} \right) \]
and inequality \eqref{Case1} yields that
\[ |{\it\Lambda}_1| < \alpha_T^{m_1-m+2}. \]
If we further assume that ${\it\Lambda}_1 > 0$, then we get
\[ 0< n_1 \left( \dfrac{\log \alpha}{\log \alpha_T} \right) -m +
\dfrac{\log((\alpha^{n-n_1}-1)/(\sqrt{5}c_\alpha))}{\log \alpha_T} < \dfrac{\alpha_T^2}{\log \alpha_T} \alpha_T^{-(m-m_1)} < 6\alpha_T^{-(m-m_1)}.\]
Again we apply Lemma \ref{Dujella} with the same $\tau$ and $M$ as in the case that ${\it \Lambda}>0$.
We use the $104$-th convergent $\frac{p}{q} = \frac{p_{104}}{q_{104}}$ of $\tau$
as before. But, in this case we choose $(A,B)=(6, \alpha_T)$ and use
\[ \mu_k = \dfrac{\log((\alpha^k-1)/(\sqrt{5}c_\alpha))}{\log \alpha_T},\]
instead of $\mu$ for each possible value of $n-n_1=k=1,2,\dots, 271$. A quick computer aid computation yields that $\varepsilon > 0.00038$ for all $1\leqslant k \leqslant 271$.
Hence, by Lemma \ref{Dujella}, we get
\[ m-m_1 < \dfrac{\log(6q/0.00038)}{\log\alpha_T} <220. \]
Thus, $n-n_1 \leqslant 271$ implies $m-m_1 \leqslant 219$.

In the case that ${\it\Lambda}_1<0$ we follow the ideas from the case that ${\it\Lambda}_1>0$. We use the same $\tau$ as in the
case that ${\it\Lambda}<0$ but instead of $\mu$ we take
\[\mu_k = \dfrac{\log(\sqrt{5}c_\alpha/(\alpha^k-1))}{\log \alpha}\]
for each possible value of $n-n_1=k=1,2,\dots, 271$.
Using Lemma \ref{Dujella} with this setting we also obtain in this case that $n-n_1 \leqslant 271$ implies $m-m_1 \leqslant 219$.

Now let us turn to the case that $m-m_1 \leqslant 212$ and let us consider inequality~\eqref{Case2}. Recall that
\[ {\it\Lambda}_2= n\log \alpha - m_1 \log \alpha_T + \log \left( \dfrac{1}{\sqrt{5}\,c_\alpha (\alpha_T^{m-m_1}-1 )} \right) \]
and let us assume that $n-n_1 \geqslant 20$. Then we have
\[ |{\it\Lambda}_2| < \dfrac{4.44 \alpha^4}{\alpha^{n-n_1}}.  \]
Assuming that ${\it\Lambda}_2 > 0$, we get
\[  0<  n \left(\dfrac{\log \alpha}{\log \alpha_T} \right)- m_1 + \dfrac{\log( 1/(\sqrt{5}\,c_\alpha (\alpha_T^{m-m_1}-1 ) ))}{\log \alpha_T} < \dfrac{4.44 \alpha^4}{(\log \alpha_T) \alpha^{n-n_1}} < 50\alpha^{-(n-n_1)} . \]
Once again we apply Lemma \ref{Dujella} with the same $\tau$ and $M$ as for the case ${\it\Lambda} > 0$ before. We take $(A,B) = (50, \alpha)$ and
\[ \mu_k=  \dfrac{\log( 1/(\sqrt{5}\,c_\alpha (\alpha_T^k-1 ) ))}{\log \alpha_T}\]
for every possible value $m-m_1=k=1,2, \dots, 212$. If we use again the $104$-th convergent of $\tau$, i.e. we put $q=q_{104}$,
then for each $k$ that yields a positive $\varepsilon$, we get $\varepsilon > 0.0012$. Therefore we get
\[ n-n_1 < \dfrac{\log(50q_{104}/0.0012)}{\log \alpha} < 280 \]
in these cases. But for $k=90$ we get a negative $\varepsilon$. In this case we consider the $105$-th convergent $\dfrac{p}{q}=\dfrac{p_{105}}{q_{105}}$ of $\tau$ instead. Let us note that
\[q_{105}= 20120013979896675119357414743592977629715414121119669783.\]
Now we obtain in the case $k=90$ that $\varepsilon > 0.46$. Thus
\[ n-n_1 < \dfrac{\log (50q_{105}/0.46)}{\log \alpha} < 275.\]

In the case that ${\it\Lambda}_2<0$ we follow again the ideas from the case that ${\it\Lambda}_2>0$. Of course we choose
\[\tau=\dfrac{\log \alpha_T}{\log \alpha}  \quad \mbox{and} \quad \mu_k=\dfrac{\log( \sqrt{5}\,c_\alpha (\alpha_T^k-1))}{\log \alpha}.\]
Applying Lemma \ref{Dujella} for all possible values of $m-m_1=k=1,\dots, 212$ also yields in this case that $n-n_1 \leqslant 279$.

Let us summarize the above computations. First we got that either $n-n_1 \leqslant 271$, or $m-m_1 \leqslant 212$.
If we assume that $n-n_1 \leqslant 271$, then we deduce that $m-m_1 \leqslant 219$,
and if we assume that $m-m_1 \leqslant 212$, then we deduce that $n-n_1 \leqslant 279$.
Altogether we obtain $n-n_1 \le 279$ and $m-m_1 \le 219$.

For the last step in our reduction process we consider inequality \eqref{Case3}. Recall that
\[ {\it\Lambda}_3 = n_1 \log \alpha - m_1 \log \alpha_T + \log\left( \dfrac{\alpha^{n-n_1}-1}{\sqrt{5}\,c_\alpha (\alpha_T^{m-m_1}-1 )} \right) . \]
Since we assume that $n \geqslant 300$, inequality \eqref{Case3} implies that
\[ |{\it\Lambda}_3| < \dfrac{3.28 \alpha^4}{\alpha^n}.  \]
Let us assume that ${\it\Lambda}_3 > 0$. Then
\[  0 < n_1 \left(\dfrac{\log \alpha}{\log \alpha_T} \right)- m_1 + \dfrac{\log\left( (\alpha^k-1)/(\sqrt{5}\,c_\alpha (\alpha_T^l-1 ))\right)}{\log \alpha_T} < \dfrac{3.28 \alpha^4}{(\log \alpha_T)\alpha^n} < 37\alpha^{-n},\]
where $(k,l) =(n-n_1, m-m_1)$. We apply Lemma \ref{Dujella} once more with the same $\tau$ and $M$ as for the case when ${\it\Lambda} > 0$. Moreover, we take $(A,B)=(37, \alpha)$, and put
\[ \mu_{k,l} = \dfrac{\log\left( (\alpha^k-1)/(\sqrt{5}\,c_\alpha (\alpha_T^l-1 ))\right)}{\log \alpha_T}\]
for $1\leqslant k \leqslant 279$ and $1\leqslant l \leqslant 219$. We consider the $104$-th convergent $\frac pq=\frac{p_{104}}{q_{104}}$.
For all pairs $(k,l)$ such that $\varepsilon$ is positive we have indeed $\varepsilon > 2.8 \times 10^{-6}$. Thus for these pairs $(k,l)$ Lemma \ref{Dujella}
yields that
\[ n \leqslant \dfrac{\log (37 q_{104}/0.0000028)}{\log \alpha} < 292.\]
For all the remaining pairs $(k,l)$ which yield a negative $\varepsilon$, we consider the $105$-th convergent $\frac pq=\frac{p_{105}}{q_{105}}$ instead.
And for all those pairs $(k,l)$ the quantity $\varepsilon$ is positive for this choice of $q$.
In particular, we have that $\varepsilon > 0.0018$ for all these cases, hence
\[ n \leqslant \dfrac{\log (37q_{105}/0.0018)}{\log \alpha} < 286.\]

In the case that ${\it\Lambda}_3<0$ the method is similar.  In particular we have to apply Lemma \ref{Dujella} with
\[\tau=\dfrac{\log \alpha_T}{\log \alpha}  \quad \mbox{and} \quad
 \mu_{k,l}= \dfrac{\log\left((\sqrt{5}\,c_\alpha (\alpha_T^l-1 ))/  (\alpha^k-1)\right)}{\log \alpha}.\]
However, we obtain in this case the slightly smaller bound  $n < 289$.

Altogether our reduction procedure yields the upper bound $n \leqslant 291$. However, this contradicts our assumption that $n \geqslant 300$.
Thus Theorem \ref{Main} is proved.

\def\cprime{$'$}


 \end{document}